\documentclass[11pt,letterpaper,reqno]{amsart}
\usepackage[T1]{fontenc}
\usepackage{lmodern,amsmath,amssymb}
\usepackage[letterpaper,textwidth=6.13in,textheight=8.48in,centering]{geometry}
\usepackage[expansion=false]{microtype}
\usepackage{enumitem}
\usepackage[hidelinks,pdfauthor={K. Srinivasa Raghava},pdftitle={Duality modulo products for finite multiple T-values}]{hyperref}
\numberwithin{equation}{section}
\newtheorem{theorem}{Theorem}[section]
\newtheorem{proposition}[theorem]{Proposition}
\newtheorem{lemma}[theorem]{Lemma}
\newtheorem{corollary}[theorem]{Corollary}
\theoremstyle{remark}
\newtheorem{remark}[theorem]{Remark}
\newcommand{\Fp}{\mathbb F_p}
\newcommand{\A}{\mathcal A}
\newcommand{\ones}[1]{\{1\}^{#1}}
\DeclareMathOperator{\dep}{dep}
\setlist[enumerate,1]{label=\textup{(\roman*)},leftmargin=*,itemsep=2pt,topsep=3pt}
\title{Duality modulo products for finite multiple \(T\)-values}
\author{K. Srinivasa Raghava}
\address{Pie Mathematics Association}
\email{srinivasaraghavak@gmail.com}
\date{}
\begin{document}
\begin{abstract}
Let \(p\) be an odd prime. We express the difference between the finite multiple \(T\)-values at an index and at the reversal of its Hoffman dual as a sum of products of two lower-weight values. Each product comes from splitting the index at a comma or a plus sign. The proof uses the two rank-one terms in the coordinate formulas for the corresponding integrals. We deduce both finite families in Cheng and Zhao's duality conjecture and a quadratic generating identity for repetitions of \((1,3)\). For \(p\ge5\) we prove an unweighted and two weighted sum formulas for the triples \(T_p(a,1,q)\) with even \(a\) and odd \(q\) in weights \(4\le w\le p+1\). They give \(T_p(1,3,1,1)=\tfrac7{32}B_{p-3}^2\) for \(p=5\) or \(p\ge11\), and express weighted Bernoulli convolutions of increasing length in terms of a fixed number of products at each weight.
\end{abstract}
\maketitle

\section{Introduction}
\subsection{Finite \texorpdfstring{\(T\)}{T}-values and duality}
Kaneko and Tsumura's multiple \(T\)-values are level-two analogues of multiple zeta values in which the summation indices alternate in parity and the smallest is odd \cite{ref11}. Fix an odd prime \(p\). For a nonempty index \(s=(s_1,\ldots,s_d)\) of positive integers, write \(\dep(s)=d\) and \(|s|=s_1+\cdots+s_d\). Its finite multiple \(T\)-value is
\begin{equation}\label{eq:1.1}
 T_p(s)=2^d\!\!\!\sum_{\substack{p>n_1>\cdots>n_d>0\\ n_j\equiv d-j+1\bmod 2}}
 \frac{1}{n_1^{s_1}\cdots n_d^{s_d}}\quad\in\Fp .
\end{equation}
Equalities involving \(T_p\) are in \(\Fp\). A rational coefficient denotes its residue whenever its denominator is prime to \(p\). For a fixed index the residues \(T_p(s)\) define a class \(T_{\A}(s)\) in \(\A=(\prod_p\Fp)/(\bigoplus_p\Fp)\). Here \(p\) ranges over odd primes. The choice of a component at \(2\) has no effect in the quotient over all primes. We write \(\ones{r}\) for \(r\) ones and \(\{1,3\}^k\) for \(k\) copies of \((1,3)\), taking either string to be empty when its exponent is zero.

For classical multiple zeta values see \cite{ref5,ref19}, and for their finite analogues see \cite{ref9,ref21}. Hoffman's duality for finite multiple harmonic sums \cite{ref6} belongs to the same setting as the depth-two evaluations \cite{ref6,ref20} and the sum formulas of \cite{ref14,ref8}. At level two, parity restrictions give the multiple \(t\)-values \cite{ref7}, the multiple mixed values \cite{ref17}, and their finite analogues \cite{ref10,ref23}. Related results include weighted sums of classical triple and alternating \(T\)-values \cite{ref2,ref18}, finite Euler sums \cite{ref16,ref22}, functional equations for finite multiple polylogarithms \cite{ref15}, and weighted sums of finite multiple mixed values \cite{ref12}. Using the normalization without \(2^d\), Zhao obtained linear shuffle relations, a sum formula, and generating sets in weights below seven for finite \(T\)-values \cite{ref22}.

Cheng and Zhao proved duality for the strings \(\ones{k}\), \((n,\ones{k})\), and \((n,2,\ones{k})\) with \(n+k\) odd in the last family. They conjectured \cite[Conjecture~4.9]{ref3} that
\begin{equation}\label{eq:1.2}
 \begin{aligned}
 T_\A(\ones{\ell},3,\ones{m})&=T_\A(\ell+1,1,m+1)
     &&(\ell\ \text{odd},\ m\ \text{even}),\\
 T_\A(\{1,3\}^k,1,1)&=T_\A(2,\{1,3\}^k) &&(k\ge1).
 \end{aligned}
\end{equation}
Theorem~\ref{thm:duality} gives both finite identities in \eqref{eq:1.2}. Corollary~\ref{cor:one}(i) proves the first for \(p\ge5\) with \(p>\ell+m+1\), and Theorem~\ref{thm:alternating} proves the second for every odd prime. The weighted sum formula \eqref{eq:1.12} gives the weight-six evaluation \eqref{eq:1.18}. The symmetric assertions in \cite[Conjecture~4.9]{ref3} concern regularized symmetric \(T\)-values and remain outside the scope of these finite arguments.

\subsection{The cut formula}
Write \(s\) as a row of \(|s|\) units, joining units within an entry by plus signs and separating entries by commas. Exchanging commas and plus signs gives the Hoffman dual \(s^\vee\) \cite{ref6}. Write \(\overleftarrow{s}=(s_d,\ldots,s_1)\) and \(s^*=\overleftarrow{s^\vee}\). Reversal commutes with this exchange, so \((s^*)^*=s\) and \((\overleftarrow{s})^*=s^\vee\). For \(k\ge2\),
\begin{equation}\label{eq:1.3}
 \begin{gathered}
 (\ones{d})^*=(d),\qquad (a,b)^*=(\ones{b-1},2,\ones{a-1}),\\
 (\ones{\ell},k,\ones{m})^*=(m+1,\ones{k-2},\ell+1).
 \end{gathered}
\end{equation}
For \(1\le i<|s|\), let \(h_i\) and \(t_i\) be the indices to the left and right of the \(i\)th separator. Splitting at a plus sign divides its entry into two positive parts.

\begin{theorem}\label{thm:duality}
For every odd prime \(p\) and every index \(s\),
\begin{equation}\label{eq:1.4}
 T_p(s)-T_p(s^*)=
 \sum_{\substack{i\ {\rm comma}\\ \dep(t_i)\ {\rm even}}}T_p(h_i^*)T_p(t_i)
 -
 \sum_{\substack{i\ {\rm plus}\\ \dep(h_i^*)\ {\rm even}}}T_p(h_i^*)T_p(t_i).
\end{equation}
\end{theorem}

For \(s=(1,3)\) the only term is the split \((1,2)\mid(1)\). This gives \(T_p(1,3)=T_p(1,1,2)-T_p(1,2)T_p(1)\), which at \(p=5\) reads \(1=3-2\cdot1\). Since \(|h_i|+|t_i|=|s|\), the classes \(T_{\A}(s)\) and \(T_{\A}(s^*)\) agree modulo the rational span of products of two lower-weight values of total weight \(|s|\). The depth of \(s^*\) is \(|s|+1-\dep(s)\). Hence, modulo these products, the weight-\(w\) values are spanned by the values of depth at most \(\lfloor(w+1)/2\rfloor\).

Two special cases used below are
\begin{equation}\label{eq:1.5}
 T_p(j)=0\quad(j\ {\rm even},\ (p-1)\nmid j),\qquad
 T_p(\ones{j})=T_p(j)\quad(1\le j\le p-1).
\end{equation}
For even depth, reversal multiplies \(T_p(s)\) by \((-1)^{|s|}\) (Lemma~\ref{result:2.6}). Applying \eqref{eq:1.4} to the reversed index therefore gives a relation with \(s^\vee\). In the families \eqref{eq:1.2}, the total correction vanishes after the appropriate reversal. For the repeated string, cancellation can take place between nonzero cuts.

\subsection{Families of indices}
\begin{theorem}\label{thm:one}
Let \(k\ge2\), \(\ell,m\ge0\), and \(\ell+m\) be odd. For every prime \(p>\ell+m+1\),
\begin{equation}\label{eq:1.6}
 \begin{aligned}
 T_p(\ones{\ell},k,\ones{m})
 &=T_p(m+1,\ones{k-2},\ell+1)\\[-2pt]
 &\quad-\sum_{\substack{2\le r<k\\r\ {\rm even}}}
 T_p(\ones{r-1},\ell+1)\,T_p(k-r,\ones{m}).
 \end{aligned}
\end{equation}
\end{theorem}
For \(k=2\) the sum is empty, giving the depth-two duality in \eqref{eq:1.14} for weights at most \(p\). For the equality in \(\A\) see \cite[Theorem~4.1 and Corollary~4.7]{ref3}. Zhao's earlier proportionality formula for indices with a single \(2\) is \cite[Proposition~3]{ref22}.

\begin{corollary}\label{cor:one}
Suppose \(p\ge5\). In \textup{(i)} and \textup{(iii)} assume \(\ell,m\ge0\), \(\ell+m\) odd, and \(p>\ell+m+1\).
\begin{enumerate}
\item \(T_p(\ones{\ell},3,\ones{m})\) equals \(T_p(\ell+1,1,m+1)\) for odd \(\ell\), and \(T_p(m+1,1,\ell+1)\) for even \(\ell\).
\item If \(a\ge2\) is even, \(q\ge1\) is odd, and \(a+q\le p\), then
\[
 T_p(a,1,q)-T_p(q,1,a)=T_p(a,1)T_p(q).
\]
\item \(T_p(\ones{\ell},4,\ones{m})+T_p(\ell+1,1,1,m+1)
 =(-1)^{\ell+1}T_p(\ell+1,1)T_p(m+1,1)\).
\end{enumerate}
\end{corollary}
For even \(\ell\), reversing the three summation indices and replacing each \(n\) by \(p-n\) interchanges the parity patterns. The weight is even, so the common value in (i) also equals the \(S\)-value at \((\ell+1,1,m+1)\), defined by exchanging odd and even in \eqref{eq:1.1} \cite{ref17,ref3}. The bound on the depth is essential in general. At \(p=5,\ell=1,m=4\) the left side has depth six and vanishes, while the only triple in \(T_5(2,1,5)\) is \((3,2,1)\), giving \(1\). In (iii) the product vanishes when \(\ell m=0\). Write \(B_j\) for the Bernoulli numbers, with \(B_1=-1/2\). For \((\ell,m)=(1,2)\) the product is \(\tfrac74 B_{p-3}T_p(1,3)\), which equals \(3\) at \(p=5\).

\begin{theorem}\label{thm:alternating}
For every odd prime \(p\) and every \(k\ge0\),
\begin{equation}\label{eq:1.7}
 T_p(\{1,3\}^k,1,1)=T_p(2,\{1,3\}^k).
\end{equation}
In \(\Fp[z]\) we also have
\begin{equation}\label{eq:1.8}
 \begin{aligned}
 &\left(1+\sum_{k\ge1}T_p(\{1,3\}^k)z^k\right)
 \left(1-\sum_{k\ge1}T_p(2,\{1,3\}^{k-1},1,1)z^k\right)\\
 &\hspace{40mm}=1-z\left(\sum_{k\ge0}T_p(\{1,3\}^k,1,1)z^k\right)^2 .
 \end{aligned}
\end{equation}
\end{theorem}
The sums terminate because depth at least \(p\) gives an empty sum. For \(p\ge5\) the coefficient of \(z^2\) gives
\begin{equation}\label{eq:1.9}
 T_p(1,3,1,3)=T_p(2,1,3,1,1)+T_p(1,3)^2.
\end{equation}

\subsection{Triple sums and Bernoulli numbers}
Define \(D_j\) by
\begin{equation}\label{eq:1.10}
 \frac{u}{\sinh u}=\sum_{j\ge0}D_j\frac{u^j}{j!},\qquad D_j=(2-2^j)B_j.
\end{equation}
Since the series is even, \(D_j=0\) for odd \(j\).

\begin{theorem}\label{thm:triple}
Let \(p\ge5\) and let \(w\) be even with \(4\le w\le p+1\). In the following sums, \(q\) runs over the odd integers \(1\le q\le w-3\) and \(a=w-1-q\). Then
\begin{align}
 \sum_q T_p(a,1,q)&=T_p(1,w-1),\label{eq:1.11}\\
 \sum_{q\ge3}(aq+q-a)T_p(a,1,q)&=\frac12\!
 \sum_{\substack{b+c=w-2\\b,c\ge2\ {\rm even}}} T_p(b,1)T_p(c,1),\label{eq:1.12}\\
 \sum_q aq\,T_p(a,1,q)&=-\frac12
 \sum_{\substack{i+j=p-w+1\\ i,j\ge0\ {\rm even}}}D_iD_j.\label{eq:1.13}
\end{align}
\end{theorem}
For each \(b\) on the right of \eqref{eq:1.12}, \(T_p(b,1)=D_{p-b-1}=(2-2^{-b})B_{p-b-1}\). Hence the number of Bernoulli products in this formula depends only on \(w\). At \(w=p+1\) each triple has the explicit value \(T_p(p-q,1,q)=-q^{-2}\), as shown in Section~5.

\begin{corollary}\label{cor:values}
Let \(p\ge5\). For \(a,b\ge1\) with \(a+b\) odd and \(a+b\le p\),
\begin{equation}\label{eq:1.14}
 T_p(a,b)=T_p(\ones{b-1},2,\ones{a-1})
 =\frac{(-1)^a}{b}\binom{a+b-1}{a}(2-2^{1-a-b})B_{p-a-b}.
\end{equation}
At weights four and five we have, with \(i,j\ge0\) even in the convolution,
\begin{align}
 \begin{gathered}
 T_p(1,3)=T_p(3,1)=T_p(2,1,1)=-\tfrac14\!\sum_{i+j=p-3}D_iD_j,\\
 T_p(2,2)=-4T_p(1,3),
 \end{gathered}\label{eq:1.15}\\
 \begin{gathered}
 T_p(1,1,2)=T_p(1,3)-\tfrac74T_p(1)B_{p-3},\\
 T_p(1,2,1)=-4T_p(1,3)+\tfrac74T_p(1)B_{p-3},
 \end{gathered}\label{eq:1.16}\\
 T_p(1,3,1)=T_p(2,1,2)+T_p(1)T_p(1,3).\label{eq:1.17}
\end{align}
At weight six, for \(p=5\) or \(p\ge11\),
\begin{align}
 \begin{aligned}
 T_p(1,3,1,1)=T_p(1,1,3,1)=T_p(2,1,3)
 &=\tfrac1{14}T_p(2,1)^2\\&=\tfrac7{32}B_{p-3}^2,
 \end{aligned}\label{eq:1.18}\\
 \begin{gathered}
 T_p(3,1,2)=-\tfrac{21}{32}B_{p-3}^2,\\
 T_p(3,1,1,1)=T_p(1,1,1,3)=T_p(4,1,1)=T_p(1,5)-\tfrac7{32}B_{p-3}^2.
 \end{gathered}\label{eq:1.19}
\end{align}
At weights eight and ten, for \(p\ge7\) and \(p\ge11\), respectively,
\begin{align}
 13T_p(2,1,5)+11T_p(4,1,3)=\tfrac{217}{64}B_{p-3}B_{p-5},\label{eq:1.20}\\
 \begin{aligned}
 19T_p(2,1,7)+21T_p(4,1,5)+15T_p(6,1,3)
 =\tfrac{889}{256}B_{p-3}B_{p-7}+\tfrac{961}{512}B_{p-5}^2 .
 \end{aligned}\label{eq:1.21}
\end{align}
\end{corollary}
The depth-two evaluation in \eqref{eq:1.14} for \(a+b<p\) is known from \cite[Theorem~4.1]{ref3} and \cite{ref23}. At \(a+b=p\), the binomial coefficient in \eqref{eq:1.14} is \(\binom{p-1}{a}=(-1)^a\) in \(\Fp\), giving the common residue \(1/b\). The relation \(T_p(2,2)=-4T_p(1,3)\) follows from the linear shuffle relations \cite{ref22,ref12}. With \eqref{eq:1.5}, formulas \eqref{eq:1.15}--\eqref{eq:1.16} express all eight values of weight four for \(p\ge7\) in terms of \(T_p(1,3)\) and \(T_p(1)B_{p-3}\). In weight six a multiple of \(B_{p-3}^2\) also occurs in Zhao's three-element generating set, together with a multiple of \(T_p(1)B_{p-5}\) \cite[Theorem~5]{ref22}. The coefficient of \(T_p(2,1,3)\) in \eqref{eq:1.12} is \(7\). At \(p=7\), Proposition~\ref{result:5.3} gives \(T_7(2,1,3)=3\), while \(\tfrac7{32}B_4^2=0\) in \(\mathbb F_7\).

\begin{corollary}\label{cor:convolution}
Let \(p\ge5\), let \(w\) be even with \(6\le w\le p+1\), and put \(n=p-w+1\). With \(a=w-1-q\), define
\[
 Q_w(j)=\sum_{\substack{3\le q\le w-3\\q\ {\rm odd}}}
 \frac{aq+q-a}{a!\,q!}\prod_{\substack{1\le s\le w-2\\s\ne q}}(j+s).
\]
Then
\begin{equation}\label{eq:1.22}
 \sum_{\substack{i+j=n\\i,j\ge0\ {\rm even}}}Q_w(j)D_iD_j
 =-\frac12\sum_{\substack{b+c=w-2\\b,c\ge2\ {\rm even}}}D_{p-b-1}D_{p-c-1}.
\end{equation}
At weights six and eight, with \(i,j\ge0\) even, the specializations are
\begin{align}
 \sum_{i+j=p-5}(j+1)(j+2)(j+4)D_iD_j
 &=-\tfrac{21}{8}B_{p-3}^2,\label{eq:1.23}\\
 \begin{split}
 \sum_{i+j=p-7}(j+1)(j+2)(j+4)(j+6)&(47j+196)D_iD_j\\
 &=-\tfrac{9765}{8}B_{p-3}B_{p-5}.
 \end{split}\label{eq:1.24}
\end{align}
Here \eqref{eq:1.23} holds for \(p=5\) or \(p\ge11\), and \eqref{eq:1.24} holds for \(p\ge7\).
\end{corollary}
The indices of the \(D\)'s in \eqref{eq:1.22} are at most \(p-w+1\) on the left and at least \(p-w+3\) on the right. Thus a convolution of increasing length reduces modulo \(p\) to a fixed number of products at each weight. For rational Bernoulli convolution identities see Miki and Gessel \cite{ref13,ref4}.

The proof of Theorem~\ref{thm:duality} uses the truncated exponential coordinates of Section~2. The two integrations become division by an index and summation over indices of the opposite parity, each with a rank-one correction. Section~3 converts these corrections into cuts. Section~4 applies them to the two families above. Section~5 uses the same coordinates as values of generating polynomials to obtain the triple sums.

\section{Truncated series and the exponential basis}\label{sec:series}

\subsection{The exponential basis}
Fix an odd prime \(p\). Let \(\mathbb Z_{(p)}\) be the ring of rational
numbers with denominator prime to \(p\), and put
\[
 R=\left\{\sum_{n\ge0}g_nu^n\in\mathbb Q[[u]]:\ 
              g_n\in\mathbb Z_{(p)}\ \text{for }0\le n<p\right\}.
\]
Let \(R_+=\{g\in R:g(0)=0\}\). For example, \(u^p/p\) lies in \(R_+\),
although its quotient by \(u\) has coefficient \(1/p\) in degree \(p-1\).
Reduction of the first \(p\) coefficients gives a surjective ring map
\(R\to V=\mathbb F_p[u]/(u^p)\), under which \(R_+\) maps onto \(uV\).
We use the same letter for a series and its reduction when the ring is
clear. The notation \([u^n]\) takes a coefficient in the ring in which
the series is written. In identities in \(\Fp\) we take the residue
of each \(p\)-integral coefficient. Let \(E\) be the image of \(\exp(u)\). Since \(E\in1+uV\),
we have \(E^p=1\). Thus \(E^k\) is the image of \(\exp(ku)\) and depends
only on \(k\bmod p\). Write \(e_k=E^k-1\).

\begin{lemma}\label{result:2.1}
The \(E^k\), \(0\le k<p\), form a basis of \(V\), and the \(e_k\),
\(1\le k<p\), form a basis of \(uV\). Moreover
\begin{gather*}
 [u^{p-1}]E^0=0,\qquad [u^{p-1}]E^k=-1\quad(1\le k<p),\\
 u=-\sum_{k=1}^{p-1}\frac{e_k}{k},\qquad
 \tanh(u/2)=-\sum_{k=1}^{p-1}(-1)^k E^k\quad\text{in }V .
\end{gather*}
\end{lemma}
\begin{proof}
Since \(E-1\) is \(u\) times a unit, the powers \((E-1)^j\), \(0\le j<p\),
form a basis of \(V\). Expanding them gives the first basis assertion,
and subtracting \(E^0\) from the other basis elements gives the second.
For \(k\ne0\bmod p\) the top coefficient is \(k^{p-1}/(p-1)!=-1\)
by Fermat's and Wilson's theorems. The coefficient of \(u^n\) in
\(-\sum_k e_k/k\), for \(1\le n<p\), is
\(-\sum_{k=1}^{p-1}k^{n-1}/n!\), which is \(1\) for \(n=1\) and
\(0\) otherwise. Finally, \(\tanh(u/2)=(\exp(u)-1)/(\exp(u)+1)\)
belongs to \(R_+\), since the denominator has constant term \(2\).
Multiplying \((E-1)/2\) by the geometric sum
\(\sum_{k=0}^{p-1}(-1)^k E^k=2/(1+E)\) gives the last formula.
\end{proof}

For \(g=\sum_{k=1}^{p-1}\varphi(k)e_k\in uV\), call \(\varphi\) its
coordinate function. In sums and parity conditions, coordinate indices are
always the integer representatives \(1,\ldots,p-1\). In particular,
\([u^{p-1}]g=-\sum_k\varphi(k)\).

\subsection{The two integration operators}
Define formal integrals by
\[
 I_1g=\int_0^u g(v)\,dv \quad(g\in R),\qquad
 I_0g=\int_0^u g(v)\frac{dv}{\sinh v}\quad(g\in R_+).
\]
For \(g=\sum_r g_ru^r\) in the respective domains and \(1\le n<p\),
\[
 [u^n]I_1g=\frac{g_{n-1}}n,\qquad
 [u^n]I_0g=\frac1n\sum_{r=1}^n g_r\frac{D_{n-r}}{(n-r)!}.
\]
Indeed, the second formula follows from \eqref{eq:1.10}. Through degree
\(p-2\), the series \(\sinh(u)/u\) has \(p\)-integral coefficients and
constant term \(1\). Its reciprocal has the same integrality through
that degree, which covers every \(D_{n-r}/(n-r)!\) above. Both integrals
therefore lie in \(R_+\), and their reductions depend only on the
reduced input coefficients. This gives maps
\(I_1:V\to uV\) and \(I_0:uV\to uV\), with \(I_1(u^{p-1})=0\)
in \(V\).

\begin{proposition}\label{result:2.2}
If \(g\in uV\) has coordinates \(\varphi\), the coordinates of \(I_1g\) and \(I_0g\) are
\begin{align}
 (I_1\varphi)(i)&=\frac{1}{i}\left(\varphi(i)+\sum_{k=1}^{p-1}\varphi(k)\right),\label{eq:2.1}\\
 (I_0\varphi)(i)&=
 \begin{cases}
 \dfrac{2}{i}\displaystyle\sum_{\substack{k>i\\k\ {\rm even}}}\varphi(k),&i\text{ odd},\\[4pt]
 -\dfrac{2}{i}\displaystyle\sum_{\substack{k<i\\k\ {\rm odd}}}\varphi(k),&i\text{ even}.
 \end{cases}\label{eq:2.2}
\end{align}
\end{proposition}
\begin{proof}
It suffices to integrate each \(e_k\) using its lift \(\exp(ku)-1\in R_+\).
The first integral gives
\begin{equation}
 I_1e_k=\frac{e_k}{k}-u=\frac{e_k}{k}+\sum_{j=1}^{p-1}\frac{e_j}{j} .\label{eq:2.3}
\end{equation}
For even \(k\), the identity
\((e^{kv}-1)/\sinh v=2\sum_{0<j<k,\ j\ {\rm odd}}e^{jv}\)
gives \(I_0e_k=2\sum_{0<j<k,\ j\ {\rm odd}}e_j/j\).
For odd \(k\), division of \(2e^v(e^{kv}-1)\) by \(e^{2v}-1\) gives
\[
 \frac{e^{kv}-1}{\sinh v}
   =2\sum_{\substack{0<j<k\\j\ {\rm even}}}e^{jv}+1+\tanh(v/2).
\]
By Lemma~\ref{result:2.1}, \(\tanh(v/2)\) has coordinates \(-(-1)^j\), whose
sum is zero. Its integral is thus \(-\sum_j(-1)^j e_j/j\) by \eqref{eq:2.3}.
Combining this with the expansion of \(u\) gives
\[
 I_0e_k=2\sum_{\substack{0<j<k\\j\ {\rm even}}}\frac{e_j}{j}
 -\sum_{j=1}^{p-1}\bigl(1+(-1)^j\bigr)\frac{e_j}{j}
 =-2\sum_{\substack{k<j<p\\j\ {\rm even}}}\frac{e_j}{j}.
\]
Reading off the \(i\)th coordinate proves \eqref{eq:2.2}.
\end{proof}

\subsection{Generating series and evaluation}
For an index \(s=(s_1,\ldots,s_d)\) define
\[
 \Phi_{s}=I_0^{s_1-1}I_1\cdots I_0^{s_d-1}I_1(1)\in R_+.
\]
In the row of units of \(s\), a plus sign contributes \(I_0\),
a comma contributes \(I_1\), and there is one final \(I_1\).
The rightmost operator acts first.
The operator \(I_1\) raises the order by one and changes parity, while
\(I_0\) preserves parity and does not lower the order. Hence
\(\Phi_{s}\) has order at least \(d\) and
\begin{equation}
 \Phi_{s}(-u)=(-1)^d\Phi_{s}(u). \label{eq:2.4}
\end{equation}

\begin{lemma}\label{result:2.3}
With \(x=\tanh(u/2)\), one has
\begin{equation}
 \Phi_{s}(u)=2^d
 \sum_{\substack{n_1>\cdots>n_d>0\\n_j\equiv d-j+1\pmod2}}
 \frac{x^{n_1}}{n_1^{s_1}\cdots n_d^{s_d}} .\label{eq:2.5}
\end{equation}
\end{lemma}
\begin{proof}
We recall the formal integral argument of \cite{ref11} and
\cite[Lemma~4.2]{ref3}.
Write \(P_{s}(x)\) for the series on the right and
\(P_{\varnothing}=1\). If \(s_1\ge2\), applying \(x\,d/dx\)
reduces \(s_1\) by one. If \(s_1=1\), summation over the outer
index after differentiation is geometric and gives
\((1-x^2)P'_{s}(x)/2=P_{(s_2,\ldots,s_d)}(x)\).
This also holds at depth one, where the geometric series starts at
\(n_1=1\). Since \(dx/du=(1-x^2)/2\) and
\(\sinh u=2x/(1-x^2)\), these recursions become the defining
integrations for \(\Phi_{s}\). Both sides vanish at \(u=0\)
for nonempty \(s\), which proves the identity by induction.
\end{proof}

\begin{theorem}\label{result:2.4}
Let \(\varphi_{s}\) be the coordinate function of the image of
\(\Phi_{s}\) in \(uV\). Then
\[
 T_p(s)=-2\sum_{k\ {\rm odd}}\varphi_{s}(k).
\]
\end{theorem}
\begin{proof}
The inverse substitution \(u(x)=\log((1+x)/(1-x))\) has
\(p\)-integral coefficients below degree \(p\) and zero constant term.
It follows that the sum
\[
 L(g)=\sum_{n=0}^{p-1}[x^n]g(u(x))\bmod p
\]
defines a linear form on \(V\). By \eqref{eq:2.5}, \(L(\Phi_{s})=T_p(s)\).
For \(0\le k<p\), summing the first \(p\) coefficients of
\(((1+x)/(1-x))^k\) gives
\[
 L(E^k)=[x^{p-1}](1+x)^k(1-x)^{-k-1}
       =[x^{p-1}](1+x)^k(1-x)^{p-1-k}=(-1)^k
\]
in \(\mathbb F_p\), since \((1-x)^p=1-x^p\).
Consequently \(L(e_k)=(-1)^k-1\), as required.
\end{proof}

\begin{corollary}\label{result:2.5}
For an index of depth \(d\),
\(\varphi_{s}(p-k)=(-1)^d\varphi_{s}(k)\).
If \(d\) is even, then \(T_p(s)=[u^{p-1}]\Phi_{s}
=-\sum_k\varphi_{s}(k)\). If \(d\) is odd, the top
coefficient is zero.
\end{corollary}
\begin{proof}
Substitution of \(-u\) for \(u\) sends \(e_k\) to \(e_{p-k}\),
so \eqref{eq:2.4} gives the first assertion. The involution \(k\mapsto p-k\)
interchanges odd and even indices. Theorem~\ref{result:2.4} and Lemma~\ref{result:2.1} now
give both top-coefficient assertions.
\end{proof}

\begin{lemma}\label{result:2.6}
If \(s\) has even depth, then
\(T_p(\overleftarrow{s})=(-1)^{|s|}T_p(s)\).
\end{lemma}
\begin{proof}
In \eqref{eq:1.1} replace \((n_1,\ldots,n_d)\) by
\((p-n_d,\ldots,p-n_1)\). For even \(d\) the required parities
are preserved. The reversed list of exponents contributes the indicated
sign.
\end{proof}

\subsection{Recovering coordinates from top coefficients}
\begin{lemma}\label{result:2.7}
For \(h\in V\) and \(0\le\ell\le p-1\),
\[
 [u^{p-1}]I_1^\ell h
 =\frac{(-1)^\ell}{\ell!}[u^{p-1}](u^\ell h).
\]
\end{lemma}
\begin{proof}
Only the coefficient of \(u^{p-1-\ell}\) in \(h\) can contribute.
After \(\ell\) integrations its multiplier is
\[
 \frac{(p-1-\ell)!}{(p-1)!}
 =\frac{(-1)^\ell}{\ell!}\quad\text{in }\Fp .
\]
This also covers the empty product at \(\ell=0\).
\end{proof}

\begin{proposition}\label{result:2.8}
If \(g\in uV\) has coordinates \(\varphi\), then for \(1\le x<p\),
\[
 \varphi(x)=\sum_{\ell=0}^{p-2}x^\ell[u^{p-1}]I_1^\ell g .
\]
\end{proposition}
\begin{proof}
The term \(\ell=p-1\) can be included, since \(u^{p-1}g=0\).
Lemma~\ref{result:2.7} identifies the sum with \([u^{p-1}](E^{-x}g)\). For
\(1\le i<p\), the coefficient of \(u^{p-1}\) in
\(E^{-x}e_i=E^{i-x}-E^{-x}\) is \(1\) when \(i=x\) and \(0\)
otherwise, by Lemma~\ref{result:2.1}. Linearity gives the result.
\end{proof}

\section{The duality formula}\label{sec:duality}

\subsection{The coordinate operators}
For a coordinate function \(\varphi\), put
\[
 (N\varphi)(i)=\frac{\varphi(i)}i,\qquad
 (K\varphi)(i)=\frac2i\sum_{\substack{k>i\\k\not\equiv i\pmod2}}\varphi(k).
\]
Write \(\lambda(\varphi)=-\sum_k\varphi(k)\),
\(\Lambda(\varphi)=-2\sum_{k\ {\rm odd}}\varphi(k)\). Set
\(\mathbf{u}(i)=-1/i\) for all \(i\), and take \(\mathbf{c}(i)=1/i\)
for even \(i\) and \(\mathbf{c}(i)=0\) for odd \(i\). Proposition~\ref{result:2.2} gives
\begin{equation}
 I_1\varphi=N\varphi+\lambda(\varphi)\mathbf{u},\qquad
 I_0\varphi=K\varphi+\Lambda(\varphi)\mathbf{c}. \label{eq:3.1}
\end{equation}
For instance, for even \(i\) the second equality subtracts the sum over
all odd \(k\) from the sum over odd \(k>i\). Each correction in \eqref{eq:3.1}
has rank one. The vector \(\mathbf{u}\) is the coordinate function of
\(u=\Phi_{(1)}\). By Theorem~\ref{result:2.4} and Corollary~\ref{result:2.5},
\begin{equation}
 \Lambda(\varphi_{s})=T_p(s),\qquad
 \lambda(\varphi_{s})=
 \begin{cases}
 T_p(s),&\operatorname{dep}(s)\text{ even},\\
 0,&\operatorname{dep}(s)\text{ odd}.
 \end{cases}\label{eq:3.2}
\end{equation}

\begin{lemma}\label{result:3.1}
For \(r\ge1\) and \(b_1,\ldots,b_r\ge0\), put
\(M=N^{b_1}KN^{b_2}K\cdots KN^{b_r}\), with \(M=N^{b_1}\) if \(r=1\), and
\(\mathbf{b}=(b_r+1,\ldots,b_1+1)\). Then
\[
 \Lambda(M\mathbf{u})=T_p(\mathbf{b}),\qquad
 \Lambda(M\mathbf{c})=
 \begin{cases}-T_p(\mathbf{b}),&r\text{ even},\\0,&r\text{ odd}.\end{cases}
\]
\end{lemma}
\begin{proof}
Expanding \(M\psi\) at \(k_1\) gives a sum over
\(k_1<\cdots<k_r<p\) of alternating parities with summand
\[
 2^{r-1}\left(\prod_{j=1}^{r-1} k_j^{-b_j-1}\right)
 k_r^{-b_r}\psi(k_r).
\]
Applying \(\Lambda\) restricts \(k_1\) to odd integers. With
\(\psi=\mathbf{u}\), reversing the chain then gives \eqref{eq:1.1} for
\(\mathbf b\). For \(\psi=\mathbf c\) a term survives only when \(k_r\)
is even, equivalently when \(r\) is even, in which case its sign changes.
\end{proof}

\subsection{Proof of Theorem~\ref{thm:duality}}
\begin{proof}
Write the operators defining \(\Phi_s\) as \(L_1\cdots L_w\), where \(w=|s|\).
The last operator sends \(1\) to \(u\), so
\(\varphi_{s}=L_1\cdots L_{w-1}\mathbf{u}\).
At weight one this product is the identity and both cut sums are empty.
For \(i<w\), write \(L_i=M_i+R_i\) using \eqref{eq:3.1}. A comma gives
\(M_i=N\), \(R_i\varphi=\lambda(\varphi)\mathbf{u}\), while a plus
gives \(M_i=K\), \(R_i\varphi=\Lambda(\varphi)\mathbf{c}\).
Telescoping yields
\begin{equation}
 L_1\cdots L_{w-1}=M_1\cdots M_{w-1}+
 \sum_{i=1}^{w-1}M_1\cdots M_{i-1}R_iL_{i+1}\cdots L_{w-1}. \label{eq:3.3}
\end{equation}
Apply \eqref{eq:3.3} to \(\mathbf{u}\) and then apply \(\Lambda\).
The positions of \(K\) in the first term become the commas of
\(s^{\vee}\), and the intervening strings of \(N\)'s count
the plus signs within its entries. Lemma~\ref{result:3.1} reverses these entries,
so the first term gives \(T_p(s^*)\).

At a cut after unit \(i\), the suffix applied to \(\mathbf{u}\) is
\(\varphi_{t_i}\). Applied to \(\mathbf{u}\), the prefix of
\(M\)'s evaluates under \(\Lambda\) to \(T_p(h_i^*)\)
by the same argument. A comma cut therefore contributes
\(T_p(h_i^*)\lambda(\varphi_{t_i})\), which by
\eqref{eq:3.2} equals \(T_p(h_i^*)T_p(t_i)\) if \(t_i\) has even depth
and is zero otherwise. At a plus sign the prefix acts on
\(\mathbf{c}\) instead. Its number of \(K\)'s is
\(\operatorname{dep}(h_i^*)-1\), so Lemma~\ref{result:3.1} gives
\(-T_p(h_i^*)T_p(t_i)\) when
\(\operatorname{dep}(h_i^*)\) is even, and zero otherwise.
These are precisely the two sums in \eqref{eq:1.4}.
\end{proof}

\subsection{First consequences}
For even \(j\), reflection \(n\mapsto p-n\) identifies the sums of
\(n^{-j}\) over the even and the odd representatives. Hence \(T_p(j)\)
is the full sum over \(\mathbb F_p^\times\), equal to zero unless
\((p-1)\mid j\).

Zhao proved even-depth vanishing for the all-ones index in \(\A\)
\cite[Proposition~7]{ref22}. The corresponding depth-one equality
and its primewise form appear in \cite[Corollary~2.6 and Remark~2.7]{ref3}.
\begin{corollary}\label{result:3.2}
For \(1\le d\le p-1\), \(T_p(\{1\}^d)=T_p(d)\). In particular
this value vanishes for even \(d\le p-3\).
\end{corollary}
\begin{proof}
Only comma cuts occur, so \eqref{eq:1.4} gives
\[
 T_p(\{1\}^{d})=T_p(d)+
 \sum_{\substack{1\le j<d\\d-j\ {\rm even}}}
   T_p(j)T_p(\{1\}^{d-j}).
\]
For even \(d\) the first factor in every product vanishes by the preceding
power sum. For odd \(d\) the second factor vanishes by induction, since
\(d-j\) is even and at most \(p-3\).
\end{proof}
At \(d=p-1\) the common value is \(-1\). For \(d\ge p\) the all-ones
sum is empty, whereas \(T_p(d)=-1\) whenever \(p-1\) divides \(d\).

We record also the depth-one evaluation in terms of the Fermat
quotient of \(2\), as in \cite[Theorem~2.5]{ref3}.
\begin{lemma}\label{result:3.3}
One has \(T_p(1)=2(2^{p-1}-1)/p\) in \(\mathbb F_p\).
\end{lemma}
\begin{proof}
For \(1\le t<p\), \(\binom pt/p=\binom{p-1}{t-1}/t\equiv
(-1)^{t-1}/t\bmod p\). Summing gives \((2^p-2)/p\).
Since \(\sum_{t=1}^{p-1}1/t=0\) in \(\mathbb F_p\), the alternating
sum is twice the sum over odd \(t\), namely \(T_p(1)\).
\end{proof}

Part (i) and the odd-\(n+k\) case of the three-case formula in (ii)
recover the dualities in \(\A\) of \cite[Theorems~4.6 and~4.8]{ref3}.
The two possible products in (ii) also describe what happens at the boundary.
\begin{corollary}\label{result:3.4}
Let \(n\ge1\) and \(k\ge0\).
\begin{enumerate}
\item If \(n\le p-1\) and \(k\le p-2\), then
\(T_p(n,\{1\}^k)=T_p(k+1,\{1\}^{n-1})\).
\item In the same range, the difference
\(T_p(n,2,\{1\}^k)-T_p(k+1,2,\{1\}^{n-1})\) is obtained by adding
\(T_p(n)T_p(k+1,1)\) if \(k\) is odd and
\(-T_p(n,1)T_p(k+1)\) if \(n\) is even.
In particular, if \(n,k\le p-3\), this difference equals
\[
 \begin{cases}
 0,&n+k\text{ odd},\\
 -T_p(n,1)T_p(k+1),&n,k\text{ even},\\
 T_p(n)T_p(k+1,1),&n,k\text{ odd}.
 \end{cases}
\]
\end{enumerate}
\end{corollary}
\begin{proof}
In (i), a plus sign can contribute only if \(h_i=(r)\)
with \(r\) even. Then \(T_p(h_i^*)=T_p(\{1\}^r)=0\), since
\(0<r<n\le p-1\). A contributing comma leaves \(j\) ones on its right
with \(j\) even and \(1\le j\le k\le p-2\). Their value is also zero. Hence \eqref{eq:1.4} gives (i).

In (ii) the same argument removes all plus cuts inside \(n\) and all
commas with only ones on their right. If \(k\) is odd, the comma
after \(n\) contributes
\[
 T_p(\{1\}^{n})T_p(2,\{1\}^{k})=T_p(n)T_p(k+1,1).
\]
If \(n\) is even, the plus inside \(2\) contributes
\[
 -T_p(2,\{1\}^{n-1})T_p(\{1\}^{k+1})=-T_p(n,1)T_p(k+1).
\]
Here (i) applies since \(k,n-1\le p-2\), and Corollary~\ref{result:3.2}
applies since \(n,k+1\le p-1\). These are all possible terms on the
stated rectangle. If \(n,k\le p-3\) and both terms occur, the even
integers \(n,k+1\) are less than \(p-1\), so both products vanish.
This proves the three cases.
\end{proof}
The distinction at the boundary occurs already at \(p=5,n=4,k=1\).
The two values in (ii) are \(2\) and \(0\), although \(n+k\) is odd.

\section{Height-one indices and repetitions of \texorpdfstring{\((1,3)\)}{(1,3)}}

\subsection{Proof of Theorem~\ref{thm:one}}
Here \(s=(\{1\}^{\ell},k,\{1\}^{m})\) has height one, meaning that exactly one entry exceeds one. Its depth $\ell+m+1$ is
even and strictly less than $p$. A comma before $k$ leaves $\{1\}^j$
on its left for some $1\le j\le\ell$, so $h_i^*=(j)$. The index on the right
has even depth precisely when $j$ is even.
For these $j$ we have $0<j<p-1$, so $T_p(j)=0$. A comma after
$k$ leaves $\{1\}^j$ on its right for some $1\le j\le m$. If $j$
is even, Corollary~\ref{result:3.2} gives $T_p(\{1\}^j)=0$.

A plus sign splits $k$ into $r$ and $k-r$, with $1\le r<k$. On the left,
$h_i=(\{1\}^{\ell},r)$ and $h_i^*=(\{1\}^{r-1},\ell+1)$ has depth $r$.
On the right $t_i=(k-r,\{1\}^m)$. Thus \eqref{eq:1.4} subtracts exactly the
products indexed by even $r$. Since
$s^{*}=(m+1,\{1\}^{k-2},\ell+1)$, this gives \eqref{eq:1.6}.

\subsection{Proof of Corollary~\ref{cor:one}}
When $k=3$ the sum has the single index $r=2$. Corollary~\ref{result:3.2} gives
$T_p(\{1\}^{m+1})=T_p(m+1)$, so \eqref{eq:1.6} becomes
\begin{equation}
 T_p(\{1\}^{\ell},3,\{1\}^{m})
 =T_p(m+1,1,\ell+1)-T_p(1,\ell+1)T_p(m+1).
 \label{eq:4.1}
\end{equation}
For part (i), suppose first that $\ell$ is even. Then $m+1$ is even. If
$m+1\le p-3$, the factor $T_p(m+1)$ vanishes. The remaining possibility
is $m=p-2$ and $\ell=0$. Here the other factor is $T_p(1,1)=T_p(2)=0$
since $p\ge5$. This proves the identity for even $\ell$. If $\ell$ is odd,
apply this case with $m$ and $\ell$ exchanged. Lemma~\ref{result:2.6} preserves the value
on the left under reversal because its weight and depth are even.

For part (ii), take $\ell=a-1$ and $m=q-1$ in \eqref{eq:4.1}. The assumption
$a+q\le p$ gives $\ell+m+1<p$. The left side is $T_p(a,1,q)$ by part (i),
and Lemma~\ref{result:2.6} gives $T_p(1,a)=-T_p(a,1)$. Hence
\[
 T_p(a,1,q)-T_p(q,1,a)=T_p(a,1)T_p(q).
\]
For part (iii), the only even $r<4$ is again $2$. Equation \eqref{eq:1.6} gives
\[
 T_p(\{1\}^{\ell},4,\{1\}^{m})
 =T_p(m+1,1,1,\ell+1)-T_p(1,\ell+1)T_p(2,\{1\}^m).
\]
The depth-four term has odd weight and acquires a minus sign under reversal.
Reversing the first factor of the product gives
$T_p(1,\ell+1)=(-1)^\ell T_p(\ell+1,1)$. Corollary~\ref{result:3.4}(i) identifies
$T_p(2,\{1\}^m)$ with $T_p(m+1,1)$. These substitutions give part (iii).

\subsection{Proof of Theorem~\ref{thm:alternating}}
Write $\{3,1\}^k$ for $k$ copies of $(3,1)$. Reversing the order of the
separators and exchanging plus signs with commas gives
\begin{equation}
\begin{aligned}
 (1,1,\{3,1\}^k)^* &= (2,\{1,3\}^k) &&(k\ge0),\\
 (1,1,\{3,1\}^{i-1},2)^* &= (\{1,3\}^i) &&(i\ge1),\\
 (\{3,1\}^i)^* &= (2,\{1,3\}^{i-1},1,1) &&(i\ge1),\\
 (\{3,1\}^i,2)^* &= (\{1,3\}^i,1,1) &&(i\ge0).
\end{aligned}
\label{eq:4.2}
\end{equation}
Put
\[
 A_k=T_p(\{1,3\}^k,1,1),\qquad B_k=T_p(2,\{1,3\}^k),
\]
and, for $k\ge1$, put
\[
 C_k=T_p(\{1,3\}^k),\qquad B'_k=T_p(2,\{1,3\}^{k-1},1,1).
\]
Set $C_0=1$. Reversal gives $A_k=T_p(1,1,\{3,1\}^k)$ and, for $k\ge1$,
$C_k=T_p(\{3,1\}^k)$, since both indices have even depth and even weight.

Apply \eqref{eq:1.4} to $s=(1,1,\{3,1\}^k)$. By \eqref{eq:4.2}, $T_p(s^*)=B_k$.
The commas with even depth on the right leave $(1,1,\{3,1\}^i)$ on the left,
$0\le i<k$, and give $B_i C_{k-i}$. A plus cut within the $i$th $3$,
$1\le i\le k$, leaves $r=1$ or $2$ units on the left of that entry.
The index on the left has weight $4i-2+r$ and depth $2i+1$, so its image
under $*$ has depth $2i-2+r$. Only $r=2$ qualifies. The index on the right is
$(1,1,\{3,1\}^{k-i})$, and the product occurs with sign $-1$ as
$-C_i A_{k-i}$. Thus
\[
 A_k=B_k+\sum_{i=0}^{k-1} B_i C_{k-i}-\sum_{i=1}^k C_i A_{k-i}.
\]
Reindexing the last sum gives
\[
 A_k-B_k=-\sum_{i=0}^{k-1}(A_i-B_i)C_{k-i}.
\]
The sum is empty when $k=0$, and induction gives $A_k=B_k$ for every
$k\ge0$. This proves \eqref{eq:1.7} for every odd prime. Individual products may
be nonzero. For example, at $p=11$ direct evaluation of \eqref{eq:1.1} gives
$A_1=B_1=9$ and $C_1=6$. At $k=2$ the terms $B_1 C_1$ and $-C_1 A_1$
are $10$ and $-10$ in $\mathbb F_{11}$, and cancel.

Next take $s=(\{3,1\}^k)$ for $k\ge1$. Now $T_p(s^*)=B'_k$.
A comma cut after $i$ complete blocks gives $B'_i C_{k-i}$ for
$1\le i<k$. These are all commas with even depth on the right. Splitting
the next $3$ leaves $(\{3,1\}^i,r)$ on the left. Its image under $*$
has depth $2i+r$. Thus $r=2$, and \eqref{eq:4.2} identifies the
product as $-A_i A_{k-1-i}$ for $0\le i<k$. We obtain
\[
 C_k=\sum_{i=1}^k B'_i C_{k-i}-\sum_{i=0}^{k-1}A_i A_{k-1-i}.
\]
This is equality of the coefficients of $z^k$ in \eqref{eq:1.8}. The constant terms
also agree. The sums define polynomials because all indices of depth at
least $p$ give empty sums.

For $p\ge5$, $A_0=T_p(1,1)=T_p(2)=0$. The coefficient of $z$ gives
$C_1=B'_1$, or $T_p(1,3)=T_p(2,1,1)$. The coefficients of $z^2$ and $z^3$
then read
\[
 C_2=B'_2+C_1^2,\qquad C_3-B'_3=C_1(C_2+B'_2)-A_1^2.
\]
The first identity is \eqref{eq:1.9}.

\section{Triple values and Bernoulli numbers}
Throughout this section $p\ge5$.

\subsection{Sums of powers over odd integers}
For $1\le a\le p-2$ define a function on $\mathbb F_p^\times$ by
\begin{equation}
 \sigma_a(x)=-\frac1a\sum_{\substack{0\le i\le p-1-a\\ i\ {\rm even}}}
 \binom{a+i-1}{i}D_i\,x^{1-a-i}.
 \label{eq:5.1}
\end{equation}

\begin{lemma}\label{result:5.1}
For $1\le a\le p-2$ and $1\le x<p$, we have
\[
 \sigma_a(x)=
 \begin{cases}
  2\displaystyle\sum_{\substack{0<y<x\\y\ {\rm odd}}}y^{-a},&x\text{ even},\\[5pt]
  2\displaystyle\sum_{\substack{0<y<x\\y\ {\rm even}}}y^{-a}+T_p(a),&x\text{ odd}.
 \end{cases}
\]
Moreover $\sigma_a(-x)=(-1)^{a+1}\sigma_a(x)$, and for $2\le a\le p-2$,
\[
 T_p(a)=\frac{2-2^{2-a}}{a}B_{p-a}.
\]
\end{lemma}

\begin{proof}
Put $L=p-a$, so $2\le L\le p-1$. For every even integer $X\ge0$ the
power-sum formula for the Bernoulli polynomial $B_L(t)$ \cite{ref1} gives
\begin{align*}
 2\sum_{\substack{0<y<X\\y\ {\rm odd}}}y^{L-1}
 &=\frac{2^L}{L}\left[B_L\!\left(\frac{X+1}{2}\right)
                         -B_L\!\left(\frac12\right)\right]\\
 &=\frac1L\sum_{i=0}^{L-1}\binom{L}{i}D_i X^{L-i}.
\end{align*}
Indeed, $2^i B_i(1/2)=(2-2^i)B_i=D_i$. Every coefficient in the last sum
is $p$-integral since $i\le p-2$. Modulo $p$, $L^{-1}=-a^{-1}$ and
$\binom{p-a}{i}=(-1)^i\binom{a+i-1}{i}$. For $p\nmid X$ we also have
$X^{L-i}=X^{1-a-i}$. Because $D_i=0$ for odd $i$, taking $X=x$ proves
the first case of the lemma.

All exponents in \eqref{eq:5.1} have the parity of $1-a$, giving the asserted sign
under $x\mapsto-x$. For odd $x$, apply the even case to $p-x$ and put
$z=p-y$. This gives
\[
 \sigma_a(x)=-2\sum_{\substack{x<z<p\\z\ {\rm even}}}z^{-a}.
\]
The full sum over $\mathbb F_p^\times$ is zero for $1\le a\le p-2$.
Hence the sum over all even $z$ is $-T_p(a)/2$, which gives the second case.

Finally, for $a\ge2$ take $X=p+1$ in the power-sum identity. The term
with $y=p$ vanishes modulo $p$, so the left side gives $T_p(a)$.
Here $L\le p-2$, and all terms through $i=L$ are $p$-integral. Translating
the Bernoulli polynomial by $1/2$ gives
$\sum_{i=0}^L\binom{L}{i}D_i=2^L B_L(1)=2^L B_L$ since $L\ge2$.
The right side therefore reduces to
\[
 \frac{2^L B_L-D_L}{L}
 =\frac{(2^{L+1}-2)B_L}{L}
 =\frac{2-2^{2-a}}{a}B_{p-a}\qquad\text{in }\mathbb F_p.
\]
\end{proof}

\subsection{Depth two and a convolution formula}
The depth-two evaluation for $a+q<p$ is known
\cite[Theorem~4.1]{ref3}, \cite{ref23}.
The following proof includes $a+q=p$ and supplies the calculation needed
for depth three.

\begin{proposition}\label{result:5.2}
Let $a\ge1$, $1\le q\le p-2$, and suppose $a+q$ is odd and $a+q\le p$.
Then
\[
 T_p(a,q)=\frac{(-1)^a}{q}\binom{a+q-1}{a}D_{p-a-q}.
\]
In particular, $T_p(b,1)=D_{p-b-1}$ for every even $b$ with $2\le b\le p-1$.
\end{proposition}

\begin{proof}
In \eqref{eq:1.1}, let $x$ be the larger index. It is even, and the inner sum over
the smaller odd index is $\sigma_q(x)/2$. Thus
\[
 T_p(a,q)=2\sum_{\substack{0<x<p\\x\ {\rm even}}}x^{-a}\sigma_q(x)
         =\sum_{x\in\mathbb F_p^\times}x^{-a}\sigma_q(x).
\]
The second equality uses $a+q$ odd, which makes $x^{-a}\sigma_q(x)$ even
under $x\mapsto-x$. Its monomials are $x^{-e}$ with $e=a+q+i-1$ and
$0<e\le a+p-2\le2p-3$, since $i\le p-1-q$. Summation over $\mathbb F_p^\times$ retains only $e=p-1$,
with multiplier $-1$. Thus $i=p-a-q$ and
\[
 T_p(a,q)=\frac1q\binom{p-a-1}{q-1}D_{p-a-q}.
\]
Modulo $p$ the binomial coefficient is
$(-1)^{q-1}\binom{a+q-1}{q-1}$, with $(-1)^{q-1}=(-1)^a$.
\end{proof}

For $2\le j\le p$ we have $D_{p-j}=(2-2^{1-j})B_{p-j}$ in $\mathbb F_p$.
Combining the proposition with the $k=2$ case of Theorem~\ref{thm:one} gives \eqref{eq:1.14}
when $b\le p-2$. At $a+b=p$, its right side is $1/b$, since
$\binom{p-1}{a}=(-1)^a$ in $\mathbb F_p$. The remaining pair $(a,b)=(1,p-1)$
has $T_p(1,p-1)=-T_p(p-1,1)=-1$ by Lemma~\ref{result:2.6} and Proposition~\ref{result:5.2}.
Theorem~\ref{thm:one} gives the dual equality here as well.

\begin{proposition}\label{result:5.3}
Let $1\le a,q\le p-2$ with $a+q$ odd. Then
\begin{equation}
\begin{aligned}
 T_p(a,1,q)={}&T_p(a)T_p(1,q)\\
 &+\frac1{aq}\sum_{\substack{i+j=p-a-q\\i,j\ge0\ {\rm even}}}
   \binom{a+i-1}{i}\binom{q+j-1}{j}D_iD_j.
\end{aligned}
\label{eq:5.2}
\end{equation}
If $a+q>p$ the sum is empty. If $a$ is even the product $T_p(a)T_p(1,q)$ vanishes.
\end{proposition}

\begin{proof}
Fix the middle index $x$, which is even. The two odd indices lie on
opposite sides of $x$. Their sums are $(T_p(a)-\sigma_a(x))/2$ above $x$
and $\sigma_q(x)/2$ below it. For every such $q$, the definition gives
$T_p(1,q)=2\sum_{0<x<p,\ x\ {\rm even}}\sigma_q(x)/x$. Hence
\[
 T_p(a,1,q)=T_p(a)T_p(1,q)
  -2\sum_{\substack{0<x<p\\x\ {\rm even}}}x^{-1}\sigma_a(x)\sigma_q(x).
\]
Again the summand is even under $x\mapsto-x$. We may replace twice the
restricted sum by the full sum. The exponents are $-e$ with
$e=a+q+i+j-1$. The bounds $i\le p-1-a$ and $j\le p-1-q$ give
$0<e\le2p-3$, so a term survives precisely when $i+j=p-a-q$.
This is impossible if $a+q>p$. The full sum at $e=p-1$ gives a factor of $-1$, which cancels the sign
in front of it and yields \eqref{eq:5.2}. If $a$ is even, then $2\le a\le p-3$
and \eqref{eq:1.5} gives $T_p(a)=0$.
\end{proof}

The convolution in \eqref{eq:5.2} is symmetric in $a$ and $q$. Interchanging them
therefore gives a second proof of Corollary~\ref{cor:one}(ii) for $a\le p-3$, using
$T_p(1,a)=-T_p(a,1)$ for even $a$.

\subsection{The sum formula}
For $\ell,m\ge0$ let $\tau(\ell,m)\in\Fp$ be the coefficient in the reduction of the series,
\[
  \tau(\ell,m)=[u^{p-1}]\Phi_{(\{1\}^{\ell},3,\{1\}^{m})}.
\]
By Corollary~\ref{result:2.5} this equals $T_p(\{1\}^{\ell},3,\{1\}^{m})$ if
$\ell+m$ is odd and vanishes if $\ell+m$ is even. It also vanishes when
$\ell+m\ge p-1$, since the order of the series is at least $\ell+m+1$.

\begin{proposition}\label{result:5.4}
For odd $k$ with $1\le k\le p-2$ one has
\[
 \sum_{\ell+m=k}T_p(\{1\}^{\ell},3,\{1\}^{m})=2T_p(1,k+2),
\]
where $\ell,m\ge0$.
\end{proposition}
\begin{proof}
The defining operators are $I_1^\ell I_0^2 I_1^{m+1}$. For $1\le x,y<p$ let
$\Gamma(x,y)$ denote the $x$th coordinate of $I_0^2e_y$. Expanding $e_y$ in
powers of $u$ and applying Proposition~\ref{result:2.8} gives
\[
 \Gamma(x,y)=\sum_{\ell,m\ge0}\tau(\ell,m)x^\ell y^{m+1},
 \qquad e_y=\sum_{m=0}^{p-2}\frac{y^{m+1}u^{m+1}}{(m+1)!}.
\]
Two applications of \eqref{eq:2.2} give
\[
 \Gamma(x,x)=-\frac4x
 \begin{cases}
 \displaystyle\sum_{0<j<x,\ j\ {\rm odd}}j^{-1},&x\ {\rm even},\\[4pt]
 \displaystyle\sum_{x<j<p,\ j\ {\rm even}}j^{-1},&x\ {\rm odd}.
 \end{cases}
\]
Indeed $I_0e_x$ is supported on the odd $j<x$ with coordinates $2/j$
when $x$ is even, and on the even $j>x$ with coordinates $-2/j$ when $x$ is odd.
On the diagonal the preceding expansion gives
$\Gamma(x,x)=\sum_{k=0}^{p-2}c_k x^{k+1}$, with
$c_k=\sum_{\ell+m=k}\tau(\ell,m)$. The exponents $1,\ldots,p-1$
are distinct modulo $p-1$, so summation on $\mathbb F_p^\times$ gives
\[
 c_k=-\sum_{x\in\mathbb F_p^\times}x^{-(k+1)}\Gamma(x,x)
     =T_p(k+2,1)+T_p(1,k+2),
\]
the two terms coming from even and odd $x$, respectively.
For odd $k$ they agree by even-depth reversal (Lemma~\ref{result:2.6}).
\end{proof}

\begin{proof}[Proof of \eqref{eq:1.11}]
Take $k=w-3$ in Proposition~\ref{result:5.4}. Even at $w=p+1$ one has $k\le p-2$,
so Corollary~\ref{cor:one}(i) applies. It identifies the summands
with $T_p(a,1,q)$, where $q$ is odd, $1\le q\le w-3$, and $a=w-1-q$.
Each occurs twice, at $\ell=a-1$ and $\ell=q-1$. Division by $2$ gives \eqref{eq:1.11}.
\end{proof}

\subsection{The unweighted convolution of the \texorpdfstring{\(D_j\)}{coefficients}}
Throughout the proof of Theorem~\ref{thm:triple} put $n=p-w+1$, write
$t_q=T_p(a,1,q)$ for $a=w-1-q$ and $q$ odd, $1\le q\le w-3$,
and define a linear functional on $\mathbb F_p[j]$ by
\[
  \mu(P)=\sum_{\substack{0\le j\le n\\j\ {\rm even}}} P(j)D_jD_{n-j}.
\]
\begin{proof}[Proof of \eqref{eq:1.13}]
For $w<p$, Proposition~\ref{result:5.3} expresses $t_q$ by the convolution in \eqref{eq:5.2}.
At $w=p+1$, Corollary~\ref{cor:one}(i) identifies $t_q$ with
$T_p(\{1\}^{a-1},3,\{1\}^{q-1})$, where $a+q=p$.
Its depth is $p-1$, so its only chain has $n_j=p-j$ and $n_a=q$.
By Fermat's and Wilson's theorems its value is
$2^{p-1}/((p-1)!\,q^2)=-q^{-2}=1/(aq)$. This is precisely the
convolution term of \eqref{eq:5.2} with $i=j=0$, including the case $a=p-1$.
Thus the same convolution represents $t_q$ throughout the stated range.

Let $\pi_c(j)=\prod_{1\le s\le w-2,\ s\ne c}(j+s)$ for $1\le c\le w-2$.
We have $i=n-j\equiv-a-q-j\pmod p$. Thus
\[
 \binom{a+i-1}{i}
   =\frac{(-1)^{a-1}}{(a-1)!}\prod_{s=q+1}^{w-2}(j+s),\qquad
 \binom{q+j-1}{j}=\frac{1}{(q-1)!}\prod_{s=1}^{q-1}(j+s)
\]
in $\mathbb F_p$. Since $a$ is even, the convolution gives
\begin{equation}
 a!\,q!\,t_q=-\mu(\pi_q). \label{eq:5.3}
\end{equation}
Since $n\equiv1-w\pmod p$, the factor $n-j+s$ equals
$-(j+w-1-s)$ in $\Fp$. The substitution $s\mapsto w-1-s$ permutes
$1,\ldots,w-2$ and sends the omitted factor $c$ to $w-1-c$.
There are $w-3$ factors, an odd number, so
$\pi_c(n-j)=-\pi_{w-1-c}(j)$.
The change $j\mapsto n-j$ gives $\mu(\pi_c)=-\mu(\pi_{w-1-c})$.
Interpolate the constant polynomial at the distinct points
$-1,\ldots,-(w-2)$ to obtain
\[
 1=\sum_{c=1}^{w-2}\frac{\pi_c(j)}{\pi_c(-c)},\qquad
 \pi_c(-c)=(-1)^{c-1}(c-1)!(w-2-c)!.
\]
All these denominators are units since $w-2\le p-1$.
Those for $c$ and $w-1-c$ are negatives of one another, so the
two terms contribute equally to $\mu(1)$. In each pair the odd member $q$
has $\pi_q(-q)=(q-1)!(a-1)!$. Hence \eqref{eq:5.3} gives
$\mu(1)=-2\sum_q aq\,t_q$, which proves \eqref{eq:1.13}.
\end{proof}

\subsection{The weighted sum formula}
For $\ell\ge0$ let
\[
 F_\ell=I_0(u^{\ell+1})\in R_+,\qquad
 F_\ell'=\frac{u^{\ell+1}}{\sinh u}.
\]
The coefficients of $F_\ell'$ through degree $p-2$ are $p$-integral.
Since $F_m(0)=0$, this suffices to define $[u^{p-1}](F_\ell'F_m)$
modulo $p$.

\begin{lemma}\label{result:5.5}
For $\ell,m\ge0$ with $\ell+m\le p-2$ one has
\[
 \tau(\ell,m)=\frac{(-1)^{\ell+1}}{(\ell+1)!(m+1)!}
                 [u^{p-1}](F_\ell'F_m).
\]
\end{lemma}
\begin{proof}
The definition of $\Phi_s$ and Lemma~\ref{result:2.7} give the first identity below.
Differentiation in $\mathbb Q[[u]]$ gives the second.
\[
 \begin{aligned}
 \tau(\ell,m)&=\frac{(-1)^\ell}{\ell!(m+1)!}
                     [u^{p-1}](u^\ell I_0F_m),\\
 (u^{\ell+1}I_0F_m)'&=(\ell+1)u^\ell I_0F_m+F_\ell'F_m.
 \end{aligned}
\]
The derivative has top coefficient $p[u^{p-\ell-1}]I_0F_m$. This vanishes
modulo $p$, since $I_0F_m\in R_+$ and $0\le p-\ell-1<p$. Comparing
coefficients proves the formula.
\end{proof}

Over $\mathbb Q$ form the series
\[
 G(X,u)=\int_0^u\frac{e^{Xv}-1}{\sinh v}\,dv
    =\sum_{\ell\ge0}F_\ell(u)\frac{X^{\ell+1}}{(\ell+1)!}.
\]
In $[u^{p-1}](\partial_uG\,\partial_XG)$ only products $F_\ell'F_m$
with $\ell+m+1\le p-1$ occur, by their orders at zero. Their relevant
coefficients are $p$-integral, and $(\ell+1)!$ and $m!$ are units.
Let $H\in\Fp[X]$ be the reduction of this polynomial. It has degree at
most $p-1$ and zero constant term. Lemma~\ref{result:5.5} gives
\begin{equation}
 H(X)=
 \sum_{\ell+m\le p-2}(-1)^{\ell+1}(m+1)\tau(\ell,m)X^{\ell+m+1}
 . \label{eq:5.4}
\end{equation}

\begin{lemma}\label{result:5.6}
For $1\le x<p$ one has, in $\mathbb F_p$,
\[
 H(x)=\tfrac12\sigma_1(x)^2
       -2\sum_{\substack{0<r<x\\r\not\equiv x\pmod 2}}\frac{\sigma_1(r)}r .
\]
\end{lemma}
\begin{proof}
First take $x$ even and write $O=\{k:0<k<x,\ k\ {\rm odd}\}$.
The geometric sum $(e^{xv}-1)/\sinh v=2\sum_{k\in O}e^{kv}$ gives
\[
 \partial_uG(x,u)=2\sum_{k\in O}e^{ku},\qquad
 \partial_XG(x,u)=\Phi_{(2)}(u)+2\sum_{j\in O}\int_0^u v e^{jv}\,dv .
\]
In the second identity we used $\Phi_{(2)}=I_0u$. Multiplication in $V$ yields
\begin{equation}
 H(x)=2\sum_{k\in O}[u^{p-1}](E^k\Phi_{(2)})
      +4\sum_{k,j\in O}[u^{p-1}]\bigl(E^k I_1(uE^j)\bigr).
 \label{eq:5.5}
\end{equation}
For $j,k\in O$, the integers $j,k$ lie strictly between $0$ and $p$,
and $k+j$ is positive, even and less than $2p$. Thus all three are units
modulo $p$. We have
$I_1(uE^j)=uE^j/j-e_j/j^2$ and $[u^{p-1}](uE^r)=1/r$ for $r\ne0$.
Also $[u^{p-1}]E^{k+j}=[u^{p-1}]E^k=-1$ by Lemma~\ref{result:2.1}.
Thus the coefficient in the double sum is $1/(j(k+j))$. Symmetrizing gives
\[
 \sum_{k,j\in O}\frac{1}{j(k+j)}
  =\tfrac12\left(\sum_{j\in O}j^{-1}\right)^2
  =\tfrac18\sigma_1(x)^2.
\]
For $k\in O$, the coefficient in the first sum in \eqref{eq:5.5} is
coordinate $i=p-k$ of $I_0u$ by the proof of Proposition~\ref{result:2.8}. As $i$ is
even and $u$ has coordinates $-1/j$, \eqref{eq:2.2} gives this coordinate as
$(2/i)\sum_{0<j<i,\ j\ {\rm odd}}1/j=\sigma_1(i)/i$.
Since $\sigma_1$ is even, this coefficient is $-\sigma_1(k)/k$.
Substitution in \eqref{eq:5.5} proves the formula for even $x$.

For odd $x$ apply this formula to $p-x$. Equation \eqref{eq:5.4} gives
$H(p-x)=H(x)$, since only odd $\ell+m$ contribute. Also $\sigma_1$ is
even, and Lemma~\ref{result:5.1} and Corollary~\ref{result:3.2} give
\[
 \sum_{\substack{0<r<p\\r\ {\rm even}}}\frac{\sigma_1(r)}r
    =\tfrac12 T_p(1,1)=\tfrac12 T_p(2)=0 \qquad(p\ge5).
\]
The substitution $r=p-j$ now gives
\[
 \sum_{\substack{0<j<p-x\\j\ {\rm odd}}}\frac{\sigma_1(j)}j
 =-\sum_{\substack{x<r<p\\r\ {\rm even}}}\frac{\sigma_1(r)}r
 =\sum_{\substack{0<r<x\\r\ {\rm even}}}\frac{\sigma_1(r)}r,
\]
which proves the formula for odd $x$.
\end{proof}

\begin{lemma}\label{result:5.7}
For even $a$ with $2\le a\le p-1$,
\[
 T_p(1,1,a)=T_p(a,1,1)-T_p(1)T_p(a,1).
\]
\end{lemma}
\begin{proof}
Take $q=1$ in Corollary~\ref{cor:one}(ii). Its hypothesis $a+q\le p$ gives the full range.
\end{proof}

\begin{proof}[Proof of \eqref{eq:1.12}]
For $w=4$ both sums are empty. Suppose $w\ge6$ and set $a_0=w-2$,
so $4\le a_0\le p-1$. By \eqref{eq:5.4}, the coefficient of $X^{a_0}$ in $H$ is
$\sum_{m=0}^{w-3}(-1)^m(m+1)\tau(w-3-m,m)$.
For each odd $q$, Corollary~\ref{cor:one}(i) gives $t_q$ twice. The term $m=q-1$
contributes $q\,t_q$, and $m=a-1$ contributes $-a\,t_q$.
The powers $X,\ldots,X^{p-1}$ have distinct exponents modulo $p-1$.
Since $H(0)=0$, summing $x^{-a_0}H(x)$ over $\mathbb F_p^\times$ extracts
the coefficient with a minus sign even when $a_0=p-1$. Therefore
\begin{equation}
 \sum_q (a-q)t_q=\sum_{x\in\mathbb F_p^\times} x^{-a_0}H(x).
 \label{eq:5.6}
\end{equation}

In the second term of Lemma~\ref{result:5.6}, expand $\sigma_1(r)$ by Lemma~\ref{result:5.1}.
With the representatives $1\le x,r,y<p$, this gives
\[
 \begin{aligned}
 &2\sum_{x=1}^{p-1}x^{-a_0}
       \sum_{\substack{0<r<x\\r\not\equiv x\pmod 2}}\frac{\sigma_1(r)}r
 \\ &\qquad=4\sum_{\substack{x>r>y\\\text{alternating parities}}}
        \frac1{x^{a_0}ry}
  +2T_p(1)\sum_{\substack{x>r\\x\ {\rm even},\ r\ {\rm odd}}}
        \frac1{x^{a_0}r}.
 \end{aligned}
\]
The triples with odd $x$ contribute $\tfrac12T_p(a_0,1,1)$.
Replacing $(x,r,y)$ by $(p-y,p-r,p-x)$ in the triples with even $x$
gives $\tfrac12 T_p(1,1,a_0)$ with positive sign because $w$ is even.
The last sum with its prefactor is $\tfrac12 T_p(1)T_p(a_0,1)$.
By Lemma~\ref{result:5.7} their total is $T_p(a_0,1,1)=t_1$. Hence
\begin{equation}
 \sum_{x\in\mathbb F_p^\times}x^{-a_0}H(x)
   =\tfrac12\sum_{x\in\mathbb F_p^\times}x^{-a_0}\sigma_1(x)^2-t_1.
 \label{eq:5.7}
\end{equation}

By \eqref{eq:5.1}, $\sigma_1(x)=-\sum_{0\le i\le p-3,\ i\ {\rm even}}D_i x^{-i}$.
In the square, the power sums retain only pairs with
$(p-1)\mid a_0+i+j$. Since $0<a_0+i+j<3(p-1)$ there are two cases.
The first gives $i+j=p-1-a_0=n$. In the second, write $i=p-1-b$,
$j=p-1-c$, where $b,c\ge2$ are even and $b+c=a_0$. Proposition~\ref{result:5.2}
identifies $D_{p-1-b}=T_p(b,1)$. Thus
\begin{equation}
 \sum_{x\in\mathbb F_p^\times}x^{-a_0}\sigma_1(x)^2=-\mu(1)-\Pi,\qquad
 \Pi=\sum_{\substack{b+c=w-2\\b,c\ge2\ {\rm even}}}T_p(b,1)T_p(c,1).
 \label{eq:5.8}
\end{equation}
Combining \eqref{eq:5.6}--\eqref{eq:5.8} with $\mu(1)=-2\sum_q aq\,t_q$ gives
$\sum_q(aq+q-a)t_q-t_1=\Pi/2$. The coefficient at $q=1$ is $1$,
which gives exactly \eqref{eq:1.12}.
\end{proof}

\subsection{Evaluations}
\begin{proof}[Proof of Corollary~\ref{cor:values}]
Section~5.2 proves \eqref{eq:1.14}, including the boundary pair $(1,p-1)$.
Its weight-three case and Lemma~\ref{result:5.1} give
\[
 T_p(2,1)=-T_p(1,2)=\tfrac74 B_{p-3},\qquad T_p(3)=\tfrac12 B_{p-3}.
\]

In weight four, Lemma~\ref{result:2.6} and Corollary~\ref{cor:one}(i) give
$T_p(3,1)=T_p(1,3)=T_p(2,1,1)$. By \eqref{eq:1.13},
$2T_p(2,1,1)=-\tfrac12\sum_{i+j=p-3}D_iD_j$.
For either integration and each input in its domain,
$uI_bg=I_1 I_bg+I_b(ug)$ because the two series in $\mathbb Q[[u]]$ have equal
derivatives and zero constant terms. Reduce this identity modulo $u^p$ and $p$, apply it twice
and use $u^2=2I_1u$ to get
\[
 u\Phi_{(3)}=\Phi_{(1,3)}+\Phi_{(2,2)}+2\Phi_{(3,1)}.
\]
Lemma~\ref{result:2.7} gives $[u^{p-1}](u\Phi_{(3)})=-[u^{p-1}]\Phi_{(1,3)}$.
By Corollary~\ref{result:2.5} the latter is $-T_p(1,3)$, whence
$T_p(2,2)=-4T_p(1,3)$. This completes \eqref{eq:1.15}.
The cut formula for $(1,3)$ and $(1,2,1)$ gives, respectively,
\[
 \begin{aligned}
 T_p(1,1,2)&=T_p(1,3)+T_p(1)T_p(1,2),\\
 T_p(1,2,1)&=T_p(2,2)+T_p(1)T_p(2,1).
 \end{aligned}
\]
These are \eqref{eq:1.16} after substituting the preceding evaluations.

For $(1,3,1)$ the dual is $(2,1,2)$. The first comma gives the product
$T_p(1)T_p(3,1)$. The only plus sign satisfying the depth condition leaves $(1,1)$ on the right, whose value
is $T_p(2)=0$. The other two cuts fail the depth conditions, proving \eqref{eq:1.17}.

For weight six take $p=5$ or $p\ge11$. Equations \eqref{eq:1.12} and \eqref{eq:1.11} give
\[
 7T_p(2,1,3)=\tfrac12 T_p(2,1)^2,\qquad
 T_p(4,1,1)+T_p(2,1,3)=T_p(1,5).
\]
We may divide by $7$, obtaining $T_p(2,1,3)=\tfrac7{32}B_{p-3}^2$.
Corollary~\ref{cor:one}(i) gives
\[
 \begin{aligned}
 T_p(1,3,1,1)=T_p(1,1,3,1)&=T_p(2,1,3),\\
 T_p(3,1,1,1)=T_p(1,1,1,3)&=T_p(4,1,1).
 \end{aligned}
\]
Corollary~\ref{cor:one}(ii) gives
$T_p(3,1,2)=T_p(2,1,3)-T_p(2,1)T_p(3)=-\tfrac{21}{32}B_{p-3}^2$.
This proves \eqref{eq:1.18} and \eqref{eq:1.19}. The remaining
weight-six triple with middle entry \(1\) follows from Lemma~\ref{result:5.7}
at \(a=4\) and \eqref{eq:1.14}. For the same primes this gives
\[
 T_p(1,1,4)=T_p(1,5)-\tfrac7{32}B_{p-3}^2
                    -\tfrac{31}{16}T_p(1)B_{p-5}.
\]

For $w=8$ take $p\ge7$, and for $w=10$ take $p\ge11$. The right sides of \eqref{eq:1.12} are respectively
$T_p(2,1)T_p(4,1)$ and $T_p(2,1)T_p(6,1)+\tfrac12T_p(4,1)^2$.
Inserting $T_p(4,1)=\tfrac{31}{16}B_{p-5}$ and
$T_p(6,1)=\tfrac{127}{64}B_{p-7}$ from \eqref{eq:1.14} gives \eqref{eq:1.20} and \eqref{eq:1.21}.
\end{proof}

\begin{proof}[Proof of Corollary~\ref{cor:convolution}]
Insert $t_q=-\mu(\pi_q)/(a!\,q!)$ from \eqref{eq:5.3} into \eqref{eq:1.12} and use
$T_p(b,1)=D_{p-b-1}$. This gives \eqref{eq:1.22}. For $w=6$ only $q=3$ occurs,
and
\[
 Q_6(j)=\tfrac7{12}(j+1)(j+2)(j+4),\qquad
 -\tfrac12D_{p-3}^2=-\tfrac{49}{32} B_{p-3}^2.
\]
Multiplication by $12/7$ gives \eqref{eq:1.23} for $p=5$ or $p\ge11$.
At $p=7$, $i+j=2$ and $D_2=-1/3$. The left side of
\eqref{eq:1.23} is $(8+72)D_2=6$ in $\mathbb F_7$, while the right side is zero.
For $w=8$ and $p\ge7$ we have
\[
 \begin{aligned}
 Q_8(j)&=\tfrac{11}{144}\pi_3(j)+\tfrac{13}{240}\pi_5(j)\\
 &=\tfrac1{360}(j+1)(j+2)(j+4)(j+6)(47j+196).
 \end{aligned}
\]
The right side of \eqref{eq:1.22} is $-D_{p-3}D_{p-5}=-\tfrac{217}{64}B_{p-3}B_{p-5}$.
Multiplication by $360$ gives \eqref{eq:1.24}.
\end{proof}

At the endpoint $w=p+1$, equation \eqref{eq:1.12} also gives
\[
 \sum_{\substack{2\le j\le p-3\\j\ {\rm even}}}D_jD_{p-1-j}=1-2T_p(1)
 =1-4\,\frac{2^{p-1}-1}{p}\qquad(p\ge5).
\]
Indeed $a=p-q$ and $t_q=-q^{-2}$, so the summand on the left of \eqref{eq:1.12} is
$1-2/q$. Its sum over odd $3\le q\le p-2$ is $\tfrac12-T_p(1)$.
On the right, $b+c=p-1$ gives $D_{p-b-1}D_{p-c-1}=D_cD_b$.

\begin{remark}\label{result:5.8}
In the following sums $i,j$ are nonnegative and even. Equation \eqref{eq:1.15} and
the weight-six case of Theorem~\ref{thm:triple} give
\[
 \begin{aligned}
 \sum_{i+j=p-3}D_iD_j&=-4T_p(1,3) &&(p\ge5),\\
 \sum_{i+j=p-5}D_iD_j&=-8T_p(1,5)-\tfrac78 B_{p-3}^2 &&(p=5\text{ or }p\ge11).
 \end{aligned}
\]
Thus the unweighted convolutions involve the depth-two values $T_p(1,3)$ and
$T_p(1,5)$. For even $8\le w\le p+1$, \eqref{eq:1.11} and \eqref{eq:1.12}
give two linear equations in $(w-2)/2\ge3$ triple values. Even when
their right sides are specified, these two equations alone leave at
least one free parameter. At \(w=8\), including \eqref{eq:1.13} gives
the three coefficient rows \((1,1,1),(0,11,13),(6,12,10)\) in the
order \(q=1,3,5\). Their determinant is \(-34\). Thus the three right
sides determine the triples for \(p\ge7\), \(p\ne17\). At \(p=17\)
the first two rows are independent and the third is six times the
first minus the second.
\end{remark}

\end{document}